\documentclass[twoside]{article}
	
\usepackage{amsmath} 
\usepackage{amssymb} 
\usepackage{latexsym}
\usepackage{theorem}

\newtheorem{theorem}{Theorem}
\newtheorem{lemma}[theorem]{Lemma}
\newtheorem{corollary}[theorem]{Corollary}

\newenvironment{proof}[1]{\smallskip \noindent {\bf #1}}{\qed\smallskip}

\def\qed{\ifhmode\unskip\nobreak\fi\ifmmode\ifinner\else\hskip5pt\fi\fi
 \hfill\hbox{\hskip5pt\vrule width4pt height6pt depth1.5pt\hskip1pt}}

\newcommand{\co}{\colon\thinspace} %% Colon with correct spacing for
\newcommand{\ov}[1]{\ensuremath{\overline{#1}}}

\newcommand{\RR}{\ensuremath{{\mathbb R}}}  % bold R
\newcommand{\R}[1]{\ensuremath{{\mathbb R}^{#1}}}   %bold superscripted R-variable
\newcommand{\ZZ}{\ensuremath{{\mathbb Z}}}	   %bold Z   
\newcommand{\NN}{\ensuremath{{\mathbb N}}}	   %bold N   
\newcommand{\pde}[2]
	{\ensuremath{\frac{\partial #1}{\partial #2}}}	%%partial derivative 
\def\d#1dt{\frac{d#1}{dt}}    %%variable ODE left side
\newcommand{\eps}{\ensuremath{\epsilon}}

\newcommand{\Lam}{\ensuremath{\Lambda}}
\newcommand{\del}{\ensuremath{\delta}}

\def\mylabel#1{\label{#1}}
\begin{document}
\title{Jacobians and branch points of real analytic open maps}
\author{Morris W. Hirsch\thanks{This research was partially supported by a
grant from the  National Science Foundation}}
\maketitle

\bigskip
\bigskip
\small
\noindent
{\bf Summary.} The main result is that the Jacobian determinant of
an analytic open map $f\co \R n\to \R n$ does not change sign.  A
corollary of the proof is that the set of branch points of $f$ has
dimension $\le n-2$.
\normalsize

\bigskip
\small



\bigskip
\noindent 
{\bf Abbreviated title:}   Real analytic open maps

\bigskip
\bigskip
%%%%%%%%%%%%%%%%%%%%%%%%%%%%%%%%%%%%%%%%%%%%%%%%%%%%%%%%%%%%
\paragraph{\large Introduction}{~}\\[2ex]
%%%%%%%%%%%%%%%%%%%%%%%%%%%%%%%%%%%%%%%%%%%%%%%%%%%%%%%%%%%%
The main object of this paper is to prove the following result:
%%%%%%%%%%%%%%%%%%%%%%%%%%%%%%%%%%%%%%%%%%%%%%%%%%%%%%%%%%%%
\begin{theorem}		\mylabel{th:a}
%%%%%%%%%%%%%%%%%%%%%%%%%%%%%%%%%%%%%%%%%%%%%%%%%%%%%%%%%%%%
The Jacobian of a real analytic open map $f\co \R n \to\R n$ does not
change sign. 
\end{theorem} 
%%%%%%%%%%%%%%%%%%%%%%%%%%%%%%%%%%%%%%%%%%%%%%%%%%%%%%%%%%%%
One of the referees kindly pointed out that the special case of
polynomial maps was proved by Gamboa and Ronga \cite{GR96}:
%%%%%%%%%%%%%%%%%%%%%%%%%%%%%%%%%%%%%%%%%%%%%%%%%%%%%%%%%%%%
\begin{theorem}[\sc Gamboa and Ronga]		\mylabel{th:gr}
%%%%%%%%%%%%%%%%%%%%%%%%%%%%%%%%%%%%%%%%%%%%%%%%%%%%%%%%%%%%
A polynomial map in $\R n$ is open if and only if point inverses are
finite and the Jacobian does not change sign.
\end{theorem}
The proof of Theorem \ref{th:a} is very similar to methods in
\cite{GR96}, which  are easily adapted to analytic maps; but as 
Theorem \ref{th:a} does not seem to be known, a direct proof may be
useful.

\smallskip
$f\co \R n \to\R n$ denotes a (real) analytic  map in Euclidean
n-space.  We always assume $f$ is {\em open}, that is, $f$ maps open
sets onto open sets.  Denote the Jacobian matrix of $f$ at $p\in\R n$
by $df_p= \left[\pde {f_i}{x_j}(p)\right]$.  The rank of $df_p$ is
called the {\em rank} of $f$ at $p$, denoted by $\mathsf{rk}_p f$; the
determinant of $df_p$ is the {\em Jacobian} of $f$ at $p$, denoted by
$Jf (p)$.  When the analytic function $Jf\co\R n\to\RR$ is everywhere
non-negative or everywhere non-positive (in a set $X$), we say $Jf$
{\em does not change sign (in $X$)}.

The following sets are defined for any $C^1$ map 
$g\co M\to N$ between $n$-manifolds (without boundary): 
\begin{itemize}
\item
the set $R_k=\{p\in M\co \mathsf{rk}_p g\le k\}$
\item
the {\em critical set}, $C=R_{n-1}$
\item
the {\em branch set}, $B=\{p\in U\co \mbox{$g$ is not a local
homeomorphism at $p$}$\} 
\end{itemize}
Note that $B\subset C$ by
the inverse function theorem.
\noindent
When $g$ is analytic, we also define: 
\begin{itemize}
\item
the {\em critical analytic hypersurface} $H\subset C$, comprising
those points having a neighborhood in $C$ that is an analytic
submanifold of dimension $n-1$  
\item the {\em constant rank analytic hypersurface} $V\subset H$,
at which $g|H$ has locally constant rank  
\end{itemize}

The following results are byproducts of the proof of Theorem
\ref{th:a}:

%%%%%%%%%%%%%%%%%%%%%%%%%%%%%%%%%%%%%%%%%%%%%%%%%%%%%%%%%%%%
\begin{theorem}		\mylabel{th:b}
%%%%%%%%%%%%%%%%%%%%%%%%%%%%%%%%%%%%%%%%%%%%%%%%%%%%%%%%%%%%
{~} 
\begin{description}
\item [(i)] the restricted map $f|V$ has rank $n-1$,

\item [(ii)] $f$ is a local homeomorphism at every point of $V$,

\item [(iii)] $B\subset R_{n-2}$, 

\item [(iv)]  $\dim R_{n-2} \le n-2$.

\end{description}
\end{theorem}
When $n=2$, conclusions (iii) and (iv) imply $B$ is a closed discrete set; 
thus in this case $f$ is {\em light}, i.e., point inverses are 0-dimensional.
From Stoilow \cite{St28}, which topologically characterizes germs
of light open surface maps, we obtain:
%%%%%%%%%%%%%%%%%%%%%%%%%%%%%%%%%%%%%%%%%%%%%%%%%%%%%%%%%%%%
\begin{corollary}		\mylabel{th:c}
%%%%%%%%%%%%%%%%%%%%%%%%%%%%%%%%%%%%%%%%%%%%%%%%%%%%%%%%%%%%
When $n=2$, the germ of $f$ at any
point is topologically equivalent to the germ at $0$ of the
complex function $z^d$ for some integer $d\ne 0$.
\end{corollary}

A key role in our proofs is played by the following result,  Theorem
1.4 of Church \cite{Ch63}:
%%%%%%%%%%%%%%%%%%%%%%%%%%%%%%%%%%%%%%%%%%%%%%%%%%%%%%%%%%%%
\begin{theorem}[Church]		\mylabel{th:church1.4}
%%%%%%%%%%%%%%%%%%%%%%%%%%%%%%%%%%%%%%%%%%%%%%%%%%%%%%%%%%%%
If $g\co \R n\to\R n$ is $C^n$ and open with rank $\ge
n-1$ at every point, then $g$ is a local homeomorphism. 
\end{theorem}

Our results are close to some of those obtained by Church for $C^n$
maps.  It is interesting to compare Theorem \ref{th:b}(iii) and
Corollary \ref{th:c} to the following results from  paragraphs 1.5 to
1.8 of  his paper \cite{Ch63}:

%%%%%%%%%%%%%%%%%%%%%%%%%%%%%%%%%%%%%%%%%%%%%%%%%%%%%%%%%%%%
\begin{theorem}	[Church]	\mylabel{th:church}
%%%%%%%%%%%%%%%%%%%%%%%%%%%%%%%%%%%%%%%%%%%%%%%%%%%%%%%%%%%%
Let $g\co M\to N$ be a $C^n$ map between n-manifolds.
\begin{description}
\item [(i)] If $M=N=\R n$ and $g$ is light,   the following
conditions are equivalent: 
\begin{description}
\item [(a)]  $g$ is open

\item [(b)]  $Jg$ does not change sign

\item [(c)] $B\subset R_{n-2}$. 
\end{description}
\item [(ii)]  If $M$ is compact and $g$ is open, then $g$ is light.
\end{description}
\end{theorem}

%%%%%%%%%%%%%%%%%%%%%%%%%%%%%%%%%%%%%%%%%%%%%%%%%%%%%%%%%%%%
\paragraph{\large Proofs}  
%%%%%%%%%%%%%%%%%%%%%%%%%%%%%%%%%%%%%%%%%%%%%%%%%%%%%%%%%%%%
%%%%%%%%%%%%%%%%%%%%%%%%%%%%%%%%%%%%%%%%%%%%%%%%%%%%%%%%%%%%
\begin{lemma}		\mylabel{th:lem}
%%%%%%%%%%%%%%%%%%%%%%%%%%%%%%%%%%%%%%%%%%%%%%%%%%%%%%%%%%%%
Assume the critical set of $f$ is $C=Jf^{-1}(0)=\R {n-1}\times
\{0\}$, and $f|C$ has constant rank $k,\;0\le k\le n-1$.  Then $f$ is
a local homeomorphism, $k=n-1$, and $Jf$ does not change sign in $\R
n$.
\end{lemma}
%%%%%%%%%%%%%%%%%%%%%%%%%%%%%%%%%%%%%%%%%%%%%%%%%%%%%%%%%%%%
\begin{proof}
%%%%%%%%%%%
It suffices to prove that the conclusion holds in some neighborhood of each
point, which we may take to be the origin.  

It is convenient to denote points of $\R n$ 
as $(y,t)\in\R {n-1}\times\RR$. 

By the rank theorem we assume that in some open cubical  neighborhood
$N$  of the origin,
%%%%%%%%%%%%%%%%%%%%%%%%%%%%%%%%%%%%%%%%%%%%%%%%%%%%%%%%%%%%
\begin{equation}		\label{eq:5}
%%%%%%%%%%%%%%%%%%%%%%%%%%%%%%%%%%%%%%%%%%%%%%%%%%%%%%%%%%%%
f_i(y,0)\equiv 0, \quad  i=k+1,\dots,n
\end{equation}
Identifying $N$ with $\R n$ by an analytic diffeomorphism, we assume
this holds for all $y\in\R n$.

Because $f$ is analytic and open, there is a dense open set
$\Lambda\subset \R {n-1}$ such that for every $y\in\Lam$, the map
$t\mapsto f_n (y,t)$ is not constant on any interval.  For each $y\in
\Lam$ there exists a maximal integer $\mu (y) \ge 0$ such that
\[  0 < j <\mu (y) \implies
\left(\frac{\partial}{\partial t}\right)^{j} f_n\: (y,0)= 0,\quad
\]   
Fix $y_*\in \Lam$ such that the function $\mu\co \Lam\to \NN$ takes
its minimum value $m$ at $y_*$.  Then $\mu= m$ in a precompact
open neighborhood $W\subset \Lam$ of $y_*$. 

By Taylor's theorem there exists  $\eps>0$ such that for $(y,t)$ in
the open set
\[N= W\times\, \,]-\eps,\eps\, [\,\,\subset\R {n-1}\times \RR\]
we have
%%
%%%%%%%%%%%%%%%%%%%%%%%%%%%%%%%%%%%%%%%%%%%%%%%%%%%%%%%%%%%%
\begin{equation}		\label{eq:W}
%%%%%%%%%%%%%%%%%%%%%%%%%%%%%%%%%%%%%%%%%%%%%%%%%%%%%%%%%%%%
f_n (y,t)=t^m H(y,t), \quad H (y,t)\ne 0
\end{equation}
%%

%\smallskip
{\em Claim:} If $k\le n-2$ and $(y_0,t_0)\in N$ is such that $f_n
(y_0,t_0)=0$, then $f_{n-1} (y_0,t_0)=0$.  For $t_0=0$ by
(\ref{eq:W}), and $k\le n-2$ implies $f_{n-1}(y_0,0)=0$ by
(\ref{eq:5}).

\smallskip
Now we assume $k\le n-2$ and reach a contradiction.  Since $f(N)$ is
open and contains
\[f(y_*,0)= (a_1,\dots,a_{n-2},0,0),\]
$f (N)$ also contains  points $(a_1\dots,a_{n-2},\del,0)$ with
$\del>0$.  But this contradicts the claim.

As $f$ has rank $n-1$ at every point of the critical set, 
$f$ must be a local homeomorphism by Theorem \ref{th:church1.4}.
Therefore for every $p$, the induced homomoprhism of homology groups
\[\ZZ=H_n (\R n,\R n \setminus\{p\})\to
H_n (\R n,\R n \setminus\{f(p)\})=\ZZ\] is an isomorphism, hence is
multiplication by a number $\del (p)\in\{+1, -1\}$.  

Homology theory implies that each of the two level sets of $\del \co\R
n \to \{+1,-1\}$ is open.  As $\R n$ is connected, $\del (p)$ is
constant.  As $\del (p)$ is the sign of $Jf (p)$ if $Jf(p)\ne 0$, we
have proved $Jf$ does not change sign.
\end{proof}
%%%%%%%%%%%%%%%%%%%%%

%%%%%%%%%%%%%%%%%%%%%
\subsubsection{Proof of Theorem \ref{th:a}}
%%%%%%%%%%%%%%%%%%%%%%%%%%%%%%%%%%% 
For any set $Y\subset \R n$, we say {\em the local theorem holds in
$Y$} if every point of $Y$ has a neighborhood in $Y$ in which $Jf$
does not change sign.

%%%%%%%%%%%%%%%%%%%%%%%%%%%%%%%%%%%%%%%%%%%%%%%%%%%%%%%%%%%%
\begin{lemma}		\mylabel{th:d}
%%%%%%%%%%%%%%%%%%%%%%%%%%%%%%%%%%%%%%%%%%%%%%%%%%%%%%%%%%%%
If the local theorem holds in a connected set $Y$, then $Jf$ does not
change sign in the closure $\ov Y$.
\end{lemma} 
%%%%%%%%%%%%
\begin{proof}
It suffices to prove $Jf$ does not change sign in $Y$, because $Jf$ is
continuous.  
Define $Y_+,\:Y_-$ to be the subsets of $Y$ where $Jf$ is
respectively $\ge 0$ and $\le 0$.  These sets are closed in $Y$ by
continuity of $Jf$, and  open in $Y$ by hypothesis.  As $Y$ is
connected, either $Y=Y_+$ or $Y=Y_-$.
\end{proof}

The local theorem obviously holds in the set $\R n\setminus C$ of
noncritical points.  By Lemma \ref{th:lem},  it also holds in
relatively open analytic hypersurface $V\subset C$ defined in the
introduction.  It remains to prove that every point of $C\setminus V$
has a neighborhood in which $Jf$ does not change sign.

%%%%%%%%%%%%%%%%%%%%%%%%%%%%%%%%%%%%%%%%%%%%%%%%%%%%%%%%%%%%
\begin{lemma}		\mylabel{th:e}
%%%%%%%%%%%%%%%%%%%%%%%%%%%%%%%%%%%%%%%%%%%%%%%%%%%%%%%%%%%%
Every point $p\in C\setminus V$ has a neighborhood $X_p \subset
C\setminus V$ that is an analytic variety of dimension $\le n-2$.
\end{lemma}
\begin{proof}
Write
\[C\setminus V= (C\setminus H)\cup (H\setminus V)\]
Suppose $p\in C\setminus H$.  In this case we take $X_p$ to be the
union of the variety $C_{\mathsf{sing}}$ of singular points of $C$ and
those connected components of $C\setminus C_{\mathsf{sing}}$ having
dimension $\le n-2$.  

Suppose $p\in H\setminus V$, or equivalently:
$p\in H$ and some minor determinant of $df$ vanishes at $p$ but not
identically in any neighborhood of $p$ in $H$.  We take $X_p$ to be
the intersection of $C$ with the union of the zero sets of such
minors.
\end{proof}

Now consider any point $p\in C\setminus V$.  By Lemma \ref{th:e}, $p$
has a neighborhood $X_p\subset C\setminus V$ that is the union of
finitely many smooth submanifolds of $\R n$ having dimensions $\le
n-2$.

Choose a connected open neighborhood $N_p\subset \R n$ of $p$ such that
$N_p\cap (C\setminus V)=N_p\cap X_p$.  Then 
$N_p\setminus X_p\; \subset\; (\R n\setminus C)\cup V$. 
Therefore the local theorem holds in $N_p\setminus X_p$.

Now $N_p\setminus X_p$ is connected, by a standard general position
argument.  Therefore from Lemma \ref{th:d}, with $Y=N_p\setminus X_p$,
we infer that $Jf$ does not change sign in  $\ov {N_p\setminus X_p}$,
which equals $\ov{N_p}$ because $X_p$ is nowhere dense.  
This completes the proof of Theorem \ref{th:a}.  

%%%%%%%%%%%%%%%%%%%%%%%%%%%%%%%%%%%%%%%%%%%%%%%%%%
\subsubsection{Proof of Theorem \ref{th:b}}
%%%%%%%%%%%%%%%%%%%%%%%%%%%%%%%%%%%%%%%%%%%%%%%%%%
Parts (i) and (ii) of Theorem 2 are proved by applying Lemma
\ref{th:lem} locally.  Lemma \ref{th:e} implies (iii), because
$B\subset C\setminus V$ by (ii).
For (iv), suppose $\dim R_{n-2}=n-1$.  Then the variety $R_{n-2}$
contains an analytic hypersurface, which must meet $V$. As
$R_{n-2}\subset C$, this implies $R_{n-2}\cap V\ne \varnothing$,
contradicting (i).  \qed

%%%%%%%%%%%%%%%%%%%%%%%%%%%%%%%%%%%%%%

\bigskip
\begin{flushleft}
Morris W. Hirsch\\[1ex]
Professor Emeritus\\
Department of Mathematics\\
University of California, Berkeley\\[1ex]
Honorary Fellow\\
Department of Mathematics\\
University of Wisconsin, Madison\\[1ex]

{\tt mwhirsch@chorus.net}

\end{flushleft}


\begin{thebibliography}{9}
%%%%%%%%%%%%%%%%%%%%%%%%%%%%%%%%%%%%%%
\bibitem{Ca99} {\sc A. Casson}, personal communication, 1997.

%% \bibitem{Whi72} {\sc H. Whitney}, {\em  Complex analytic varieties},  New
%%York: Addison Wesley (1972). 

\bibitem{Ch63} {\sc P. Church}, {\em Differentiable open maps on
manifolds},  Trans. Amer. Math. Soc.  {\em 109} (1963), 87--100 

\bibitem{GR96} {\sc J. Gamboa \& F. Ronga}, {\em On open real
polynomial maps}, J. Pure Appl. Algebra  {\bf 110} (1996), 297--304

\bibitem{St28} {\sc S. Stoilow}, {\em Sur les transformations continues et
la topologie des fonctions analytiques}, Ann. Sci. \'{E}cole
Norm. Sup. III  {\em 45} (1928), 347--382

\end{thebibliography}
\end{document}